\documentclass[a4paper, 10pt]{article}
\usepackage{amsmath,amsthm,amssymb,amsfonts,amscd, graphicx}
\usepackage{graphicx}
\usepackage{tikz-cd}
\usepackage{caption}
\usepackage{subcaption}

\usepackage{algorithm}
\usepackage{algpseudocode}

\usepackage[all]{xy}
\usepackage[colorinlistoftodos]{todonotes}
\usepackage{float}
\usepackage{comment}
\usepackage{booktabs}
\usepackage{mathtools}

\usepackage[left=1.1 in,top=1.2 in, right=1.1in, bottom=1.2 in]{geometry}
\usepackage[backref=page]{hyperref}

\renewcommand*{\backref}[1]{}
\renewcommand*{\backrefalt}[4]{%
  \ifcase #1 %
    \relax
  \else
    [#2]%
  \fi
}

\newtheorem{theorem}{Theorem}[section]
\newtheorem{lemma}[theorem]{Lemma}

\newtheorem{corollary}[theorem]{Corollary}

\theoremstyle{definition}
\newtheorem{definition}[theorem]{Definition}
\newtheorem{remark}[theorem]{Remark}

\numberwithin{equation}{section}
\usepackage{lineno}

\title{\bfseries{Large values in time series and additive combinatorics}}

\author{A. Iosevich\thanks{Supported in part by the National Science Foundation under grant NSF DMS-2154232} \hspace{1pt} and V. Gupta}

\date{\today}

\begin{document}

\maketitle

\begin{abstract} 
It is well-known in industrial data science that large values of real-life time series tend to be structured and often follow concrete and visible patterns. In this paper, we use ideas from additive combinatorics and discrete Fourier analysis to give this heuristic a mathematical foundation. Our main tool is the Fourier ratio, a complexity measure previously used in compressed sensing, combined with a generalized version of Chang's lemma from additive combinatorics. Together, these yield a precise prediction: when the Fourier ratio of a time series is small, the set of its largest values can be additively generated by a very small set using only $\{-1,0,1\}$ coefficients. We test this prediction on US inflation data and Delhi climate data, both in their original form and after mean-centering. The numerical results confirm the predicted structure: a generating set of size $4$--$7$ suffices to span large spectra containing dozens of points, even when the Fourier ratio is large enough that our theoretical bounds become loose. These findings provide a rigorous explanation for why extreme values in real-world data are information-rich and structurally significant.
\end{abstract} 

\section{Introduction}
\label{section:introduction} 

The observation that the largest values of a time series often exhibit more 
structure than average or typical values is a recurring theme in applied 
fields, particularly in retail analytics and demand forecasting. In these 
domains, the practical challenges posed by highly variable sales data have 
led to the development of specialized methods that implicitly acknowledge 
this phenomenon. For instance, Roque et al.~\cite{Roque2019} demonstrate 
that in hierarchical retail sales data, granular, disaggregated series—which 
contain the extreme values—hold valuable information that improves forecasts 
for aggregated totals, suggesting that the fine-scale structure is 
concentrated in the tails. Similarly, Damato et al.~\cite{Damato2025} 
explicitly model intermittent, heavy-tailed demand using a Tweedie 
likelihood, recognizing that standard distributional assumptions fail to 
capture the complexity inherent in the largest, most volatile sales events. 
This methodological shift toward heavy-tailed models is further reinforced 
by Norazizi and Syafrina~\cite{Norazizi2021}, who apply Extreme Value 
Theory to retail sales data, showing that the tail of the distribution 
follows a Generalized Pareto Distribution—a distinct statistical structure 
that does not apply to the bulk of the data. Collectively, this body of 
applied work indicates that the largest values are not simply outliers to 
be smoothed or discarded, but rather carry a structured signal that is 
critical for accurate modeling and prediction.

Beyond the development of bespoke statistical methods, the recognition of 
structured tails has also manifested in practical preprocessing techniques. 
Industry practitioners routinely employ transformations such as the 
Box–Cox transform to stabilize variance in highly volatile series, a step 
that acknowledges that the raw data's structure is too complex for models 
assuming homogeneity~\cite{Imperia2025}. Such techniques aim to make the 
underlying patterns—particularly those associated with large values—more 
accessible to standard forecasting algorithms. 

The existing literature, including foundational works on compressive sensing~\cite{Candes2006, 
Foucart2013} and complexity measures such as Kolmogorov 
complexity~\cite{Li2019, Vereshchagin2010}, has largely focused on signal 
recovery under sparsity or the information content of entire sequences. 
However, a unified theoretical framework that explains \emph{why} the 
extreme values of a time series should inherently carry more structure 
than the average, and how this structure can be quantified independently 
of domain-specific assumptions, is still lacking. 

This paper provides such a framework by connecting the heuristic of 
``structured tails'' to additive combinatorics. Specifically, we show that 
when a signal has a small Fourier ratio---a quantity that measures the 
flatness of its Fourier spectrum---the set of its largest values exhibits 
a strong additive structure: there exists a small set $\Lambda$ such that 
every large value index can be expressed as a $\{-1,0,1\}$-linear combination 
of elements of $\Lambda$ (Corollary~\ref{cor:mama}). This is a precise, 
quantifiable notion of structure that goes beyond mere statistical rarity.

\subsection{What the Theory Predicts: A Preview}
\label{subsec:preview}

Before developing the technical machinery, we state the main theoretical 
prediction in plain terms to guide the reader.

\begin{quote}
\emph{If a time series $f$ has a small Fourier ratio $\mathrm{FR}(f) = \|\hat{f}\|_1 / \|\hat{f}\|_2$, then the set of indices where $|f|$ is large can be generated by a very small set $\Lambda$ using only addition, subtraction, and the coefficients $\{-1,0,1\}$. The size of $\Lambda$ is bounded by $C\eta^{-2}\mathrm{FR}(f)^2\log(\mathrm{FR}(f)^{-2}N)$, where $\eta$ controls what counts as ``large.''}
\end{quote}

This prediction has four concrete consequences that we will test numerically:

\begin{enumerate}
\item \textbf{Small FR implies small $\Lambda$:} A small Fourier ratio should imply a small generating set $\Lambda$.

\item \textbf{Larger FR regime:} When the Fourier ratio is larger (e.g., after mean-centering), the theoretical bound becomes looser, though the additive structure may persist.

\item \textbf{$\eta$ scaling:} As the threshold $\eta$ increases (making ``large'' more restrictive), the required $|\Lambda|$ should decrease as $\eta^{-2}$.

\item \textbf{Additive generation:} Every index deemed ``large'' should be expressible as a $\{-1,0,1\}$-linear combination of the elements of $\Lambda$.
\end{enumerate}

The following subsections develop the mathematical formalism that justifies 
these predictions, after which Section~\ref{sec:numerical} puts them to the 
test on real data.

\subsection{Additive combinatorial tools}
\label{subsec:tools}

\begin{definition}
    Let $f:\mathbb{Z}_N\to \mathbb{C}$ be a signal. Its Fourier transform $\hat{f}:\mathbb{Z}_N\to \mathbb{C}$ is defined as
    \begin{equation}
        \hat{f}(m) = \frac{1}{\sqrt{N}}\sum_{x\in\mathbb{Z}_N}\chi(-xm)f(x),
    \end{equation}
    where $\chi(t)= e^{\frac{2\pi i t}{N}}.$
\end{definition}

One of the key quantities in this paper, introduced in \cite{Aldahleh2025} (see also \cite{BIN2026} for a more general setting; and \cite{BandeiraLewisMixon2015} for the numerical sparsity perspective), is the Fourier Ratio, the key measure of complexity, defined as follows. 

\begin{definition} Let $f: {\mathbb Z}_N \to {\mathbb C}$. Define 
$$ \mathrm{FR}(f)=\frac{{||\widehat{f}||}_1}{{||\widehat{f}||}_2}, $$ where here, and throughout, 
$$ {||f||}_p={\left( \sum_{x \in {\mathbb Z}_N} {|f(x)|}^p \right)}^{\frac{1}{p}}.$$
\end{definition} 

\begin{definition}
    Let $f:\mathbb{Z}_N\to \mathbb{C}$ and fix $\eta >0$. The large spectrum is defined as
    \begin{equation}
     \Gamma := \{m\in \mathbb{Z}_N \; : \; |\hat{f}(m)|\geq \eta\|f\|_{L^2(\mu)}\},  
    \end{equation}
    where $||f||_{L^2(\mu)}=\frac{||f||_2}{\sqrt{N}}$.
\end{definition}

\begin{lemma}[Chang's Lemma (\cite{Chang02})]
    Let $A\subset \mathbb{Z}_N$ have a density $\alpha = \frac{|A|}{N},$ and for $\eta >0$, define the large spectrum set 
    $$\Gamma = \{m\in \mathbb{Z}_N \; : \; |\widehat{1_A}(m)|\geq \eta\alpha N^{\frac{1}{2}}\}.$$
    Then there exists a constant $C$ and $\Lambda\subseteq \Gamma$ with $|\Lambda|\leq C\eta^{-2}\log(\frac{1}{\alpha})$ such that every $m\in \Gamma $ is a $\{-1,0,1\}$-linear combination of elements of $\Lambda$.
\end{lemma}

\begin{theorem}[Generalized Chang's Lemma (\cite{Aldahleh2025})]\label{thm: generalizedchang}
    Let $f:\mathbb{Z}_N\to\mathbb{C}$, and for $\eta>0$, define the large spectrum set
    $$\Gamma=\left\{ m\in\mathbb Z_N \,:\, |\widehat{f}(m)| \geq \eta\|f\|_{L^2(\mu)}\right\}.$$
    Then there exists a constant $C$ and some $\Lambda\subset\Gamma$ such that every $m\in \Gamma$ is a $\{-1,0,1\}$-linear combination of elements of $\Lambda$, and moreover, 
    \begin{align}\label{eq:generalizedchang_lognorm}
        |\Lambda|&\leq C\eta^{-2}\left(\frac{\|f\|_\frac{\log N}{\log N-1}}{\|f\|_2}\right)^2\log N.
    \end{align}
    If, additionally, $\frac{\|f \|_1}{\|f\|_2} \leq \frac{1}{e}\sqrt N$, then we have
    \begin{align}\label{eq:generalizedchang_l2l1}
        |\Lambda|&\leq C\eta^{-2}\frac{\|f\|_1^2}{\|f\|_2^2}\log\left(\left(\frac{\|f\|_2}{\|f\|_1}\right)^2N\right).
    \end{align}
\end{theorem}

\begin{corollary} \label{cor:mama}
    Applying Theorem \ref{thm: generalizedchang} to $\widehat{f}$, whenever $\mathrm{FR}(f) \leq \frac{1}{e}\sqrt N$ we obtain
    \begin{align}\label{eq:simple_upper_bound}
    |\Lambda|\leq C\eta^{-2}\mathrm{FR}(f)^2\log\left(\mathrm{FR}(f)^{-2}N\right),
    \end{align}
    from \eqref{eq:generalizedchang_l2l1}, where $\Lambda$ is a set such that every
    \begin{align*}
        x\in\Gamma:=\left\{x\in\mathbb{Z}_N:|f(x)|\geq\eta\|f\|_{L^2(\mu)}\right\}
    \end{align*} 
    is a $\{-1,0,1\}$-linear combination of elements of $\Lambda$. Since the $\text{FR}(f)$ term outside the $\log$ dominates, this implies that a small $\text{FR}(f)$ means a small $\Lambda$, pointing towards additive structure in the large valued set of $f$. Alternatively, for large $N$, \eqref{eq:generalizedchang_lognorm} approaches
    \begin{align}\label{eq:general_upper_bound}
        |\Lambda|\leq C'\eta^{-2}\mathrm{FR}(f)^2\log N,
    \end{align} 
    where $C'= Ce^{-2}$, giving the same result.
\end{corollary}

In short, Corollary \ref{cor:mama} provides the mathematical foundation for our investigation: large values of a well-behaved time series are structured and lend themselves to a natural compression mechanism. 

\subsection{Organization of the Paper}
\label{subsec:organization}

The remainder of the paper is organized as follows. Section~\ref{sec:numerical} formulates testable predictions from the theory (building on the preview in Section~\ref{subsec:preview}) and validates them on US inflation data and Delhi Climate data, both in its original form and after mean-centering. Section~\ref{sec:conclusion} concludes with a discussion of implications for time series analysis and directions for future work.

\section{Numerical Experiments}
\label{sec:numerical}

In this section, we test the predictions from Section~\ref{subsec:preview} using the theoretical machinery developed in Section~\ref{subsec:tools}. Our goal is to verify Corollary~\ref{cor:mama} on real-life time series data: if the Fourier ratio $\mathrm{FR}(f)$ is small, then the set of large values $\Gamma$ should be additively structured, meaning it can be generated by a small set $\Lambda$ using only $\{-1,0,1\}$ coefficients.

\subsection{Theoretical Predictions Restated}
\label{subsec:predictions_restated}

For convenience, we restate the four predictions from Section~\ref{subsec:preview} in more precise terms:

\begin{enumerate}
\item \textbf{Prediction 1 (Small FR implies small $\Lambda$):} 
Corollary~\ref{cor:mama} states that if $\mathrm{FR}(f) \leq \frac{1}{e}\sqrt{N}$, 
then $|\Lambda| \leq C\eta^{-2}\mathrm{FR}(f)^2 \log(\mathrm{FR}(f)^{-2}N)$. 
Since $\mathrm{FR}(f)$ multiplies the entire expression, a signal with a small 
Fourier ratio should require a small generating set $\Lambda$ to span its 
large values.

\item \textbf{Prediction 2 (Larger FR leads to looser bound):} 
When $\mathrm{FR}(f)$ is larger (e.g., after mean-centering), the theoretical 
upper bound increases. This does \emph{not} necessarily mean $|\Lambda|$ will 
be larger in practice, but rather that the theory guarantees less.

\item \textbf{Prediction 3 (Threshold scaling):} 
The bound in Corollary~\ref{cor:mama} scales as $\eta^{-2}$, so larger $\eta$ 
(smaller $\Gamma$) should allow proportionally smaller $\Lambda$.

\item \textbf{Prediction 4 (Additive combination structure):} 
Every element of $\Gamma$ can be written as a $\{-1,0,1\}$-linear combination 
of elements of $\Lambda$.
\end{enumerate}

\subsection{From Theory to Practice: A Dictionary}
\label{subsec:dictionary}

To avoid confusion, we clarify how the abstract theoretical objects map to 
concrete computational quantities:

\begin{table}[H]
\centering
\caption{Mapping between theoretical constructs and empirical measurements.}
\label{tab:dictionary}
\begin{tabular}{ll}
\toprule
\textbf{Theoretical object} & \textbf{Empirical counterpart} \\
\midrule
Signal $f: \mathbb{Z}_N \to \mathbb{C}$ & Time series of length $N$ (inflation rates) \\
Large spectrum $\Gamma$ & Indices where time series value exceeds threshold \\
Generating set $\Lambda$ & Set of indices selected by greedy algorithm \\
$\{-1,0,1\}$-linear combination & Integer combinations with coefficients $-1, 0, 1$ modulo $N$ \\
Fourier ratio $\mathrm{FR}(f)$ & Computed via FFT and norms \\
Bound from Corollary~\ref{cor:mama} & Theoretical ceiling on $|\Lambda|$ (up to constant $C$) \\
\bottomrule
\end{tabular}
\end{table}

\subsection{Algorithm}
\label{subsec:algorithm}

Corollary~\ref{cor:mama} guarantees the existence of such a $\Lambda$ but does not provide a construction. To test whether this structural guarantee holds on empirical data, we implement a greedy algorithm that builds $\Lambda$ incrementally. Algorithm~\ref{alg:greedy} describes the procedure.

\begin{algorithm}[H]
\caption{Greedy Construction of a Spanning Set for the Large Spectrum}
\label{alg:greedy}
\begin{algorithmic}[1]
\Require Signal $f : \mathbb{Z}_N \to \mathbb{C}$, threshold parameter $\eta > 0$
\Ensure A set $\Lambda \subseteq \Gamma$ such that every $\gamma \in \Gamma$ is a $\{-1, 0, 1\}$-linear combination of elements of $\Lambda$ modulo $N$


\State \textbf{Step 1: Identify the large spectrum.}
\[
\Gamma = \left\{ x \in \mathbb{Z}_N : |{f}(x)| \geq \eta \cdot \frac{\|f\|_2}{\sqrt{N}} \right\}.
\]

\State \textbf{Step 2: Sort $\Gamma$ in decreasing order of $|{f}|$.}
\[
\Gamma = \{\gamma_1, \gamma_2, \ldots, \gamma_{|\Gamma|}\}, \quad |{f}(\gamma_1)| \geq |{f}(\gamma_2)| \geq \cdots \geq |{f}(\gamma_{|\Gamma|})|.
\]

\State \textbf{Step 3: Greedy span construction.}
\State Initialize $\Lambda \leftarrow \emptyset$ and $S \leftarrow \{0\}$
\For{$k = 1, 2, \ldots, |\Gamma|$}
    \If{$|S| = N$}
        \State \textbf{break} \Comment{$S = \mathbb{Z}_N$, all remaining elements trivially spanned}
    \EndIf
    \If{$\gamma_k \bmod N \in S$}
        \State skip
    \Else
        \State $\Lambda \leftarrow \Lambda \cup \{\gamma_k\}$
        \State $S \leftarrow S \cup \{(s + \gamma_k) \bmod N : s \in S\} \cup \{(s - \gamma_k) \bmod N : s \in S\}$
    \EndIf
\EndFor

\State \textbf{Step 4: Verify.}
\For{each $\gamma \in \Gamma$}
    \State Check $\gamma \bmod N \in S$
\EndFor

\State \Return $\Lambda$
\end{algorithmic}
\end{algorithm}

\begin{remark}
\label{rem:constant}
The greedy algorithm is not guaranteed to produce the minimal possible $\Lambda$, but it provides an upper bound on the size needed to span $\Gamma$. After constructing $\Lambda$, we compare its size to the theoretical bound from Corollary~\ref{cor:mama}. Because the constant $C$ in Theorem~\ref{thm: generalizedchang} is not explicit, we report the quantity $\mathrm{Bound}/C$ (i.e., the bound divided by the unknown $C$). A successful test occurs when the observed $|\Lambda|$ is smaller than this quantity, which is consistent with the existence of some absolute constant $C$. The bound over the unknown $C$ can be computed using the equation (\ref{eq:simple_upper_bound}). 
\end{remark}

\subsection{US Inflation Rate (January 2023 -- February 2025)}
\label{subsec:inflation}

\begin{figure}[H]
\centering
\includegraphics[width=1\linewidth]{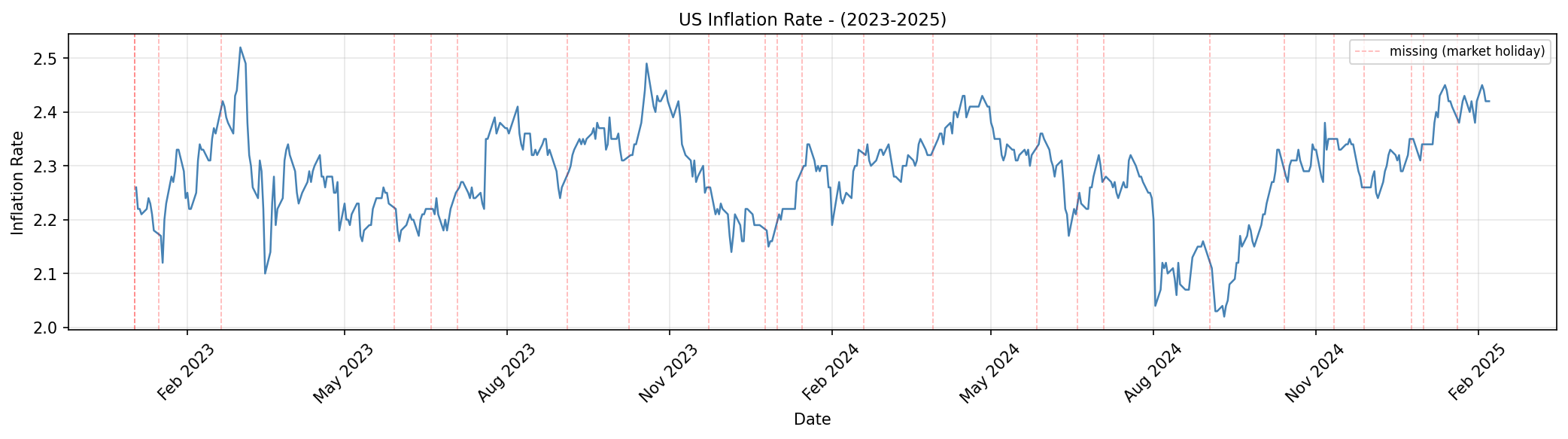}
\caption{US inflation rate from January 2023 to February 2025.}
\label{fig:inflation_rate}
\end{figure}

We first analyze the monthly US inflation rate over $N = 526$ days (Figure~\ref{fig:inflation_rate}). The signal $f$ is strictly positive. Its Fourier ratio is 
\[
\mathrm{FR}(f) = \frac{\|\widehat{f}\|_1}{\|\widehat{f}\|_2} = 1.4136,
\]
which is close to the theoretical minimum of $1$ and much smaller than 
$\frac{1}{e}\sqrt{526} \approx 8.43$. Thus, the condition for the stronger 
bound \eqref{eq:simple_upper_bound} is satisfied.

\subsubsection{Testing Prediction 1: Small FR Implies Small $\Lambda$}

Prediction~1 anticipates a compact generating set $\Lambda$. Table~\ref{tab:inflation} 
shows that $|\Lambda|$ ranges from $4$ to $7$ across different $\eta$ values, 
confirming the prediction. Moreover, even when $|\Gamma| = 75$ (for $\eta = 1.04$), 
only $|\Lambda| = 7$ basis elements are needed to span the entire set.

\subsubsection{Testing Prediction 3: $\eta$ Scaling}

The theoretical bound scales as $\eta^{-2}$. As $\eta$ increases from $1.04$ to 
$1.08$, we observe that $|\Lambda|$ decreases from $7$ to $4$, consistent with 
the predicted trend. The bound itself also decreases monotonically, as expected.

\subsubsection{Testing Prediction 4: Additive Combination Structure}

Table~\ref{tab:inflation} shows that the greedy algorithm successfully spans 
the entire $\Gamma$ in all cases (the ``Spanned'' column shows full coverage). 
This verifies Prediction~4: every element of $\Gamma$ can be expressed as a 
$\{-1,0,1\}$-linear combination of elements of $\Lambda$.

\begin{figure}[H]
\centering
\includegraphics[width=1\linewidth]{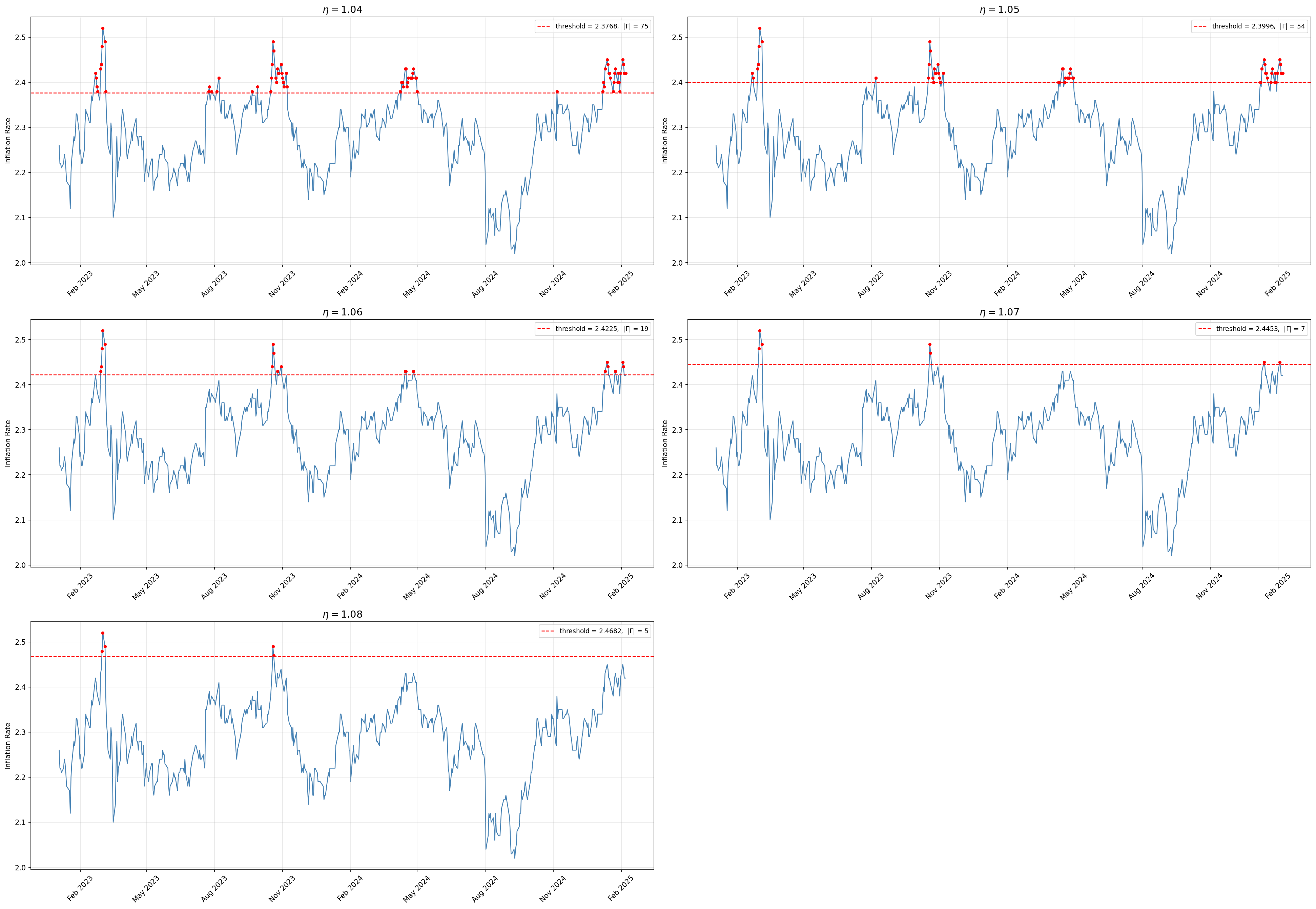}
\caption{Large spectrum $\Gamma$ of the inflation rate time series for $\eta \in \{1.04, 1.05, 1.06, 1.07, 1.08\}$. Red dots indicate points in $\Gamma$. As $\eta$ increases, $\Gamma$ condenses into tight clusters, revealing additive structure.}
\label{fig:gamma_inflation}
\end{figure}

Figure~\ref{fig:gamma_inflation} shows the locations of $\Gamma$ as red dots. 
The clustering is the visual signature of additive structure. For $\eta = 1.06$, 
$\Gamma$ contains indices around $40$--$42$, $199$--$201$, $321$--$322$, and 
$507$--$522$. These clusters are related by translations: the cluster around 
$199$--$201$ is approximately $40 + 159$ modulo $N$, and $322 \approx 40 + 282$. 
This is exactly the kind of relation that a small $\Lambda$ can capture.

\begin{table}[H]
\centering
\caption{Greedy spanning set $\Lambda$ for the large spectrum $\Gamma$ of the inflation rate time series.}
\label{tab:inflation}
\begin{tabular}{cccc}
\toprule
$\eta$ & $|\Gamma|$ & $|\Lambda|$ & Spanned \\
\midrule
$1.04$ & $75$ & $7$ & $\checkmark$ \\
$1.05$ & $54$ & $7$ & $\checkmark$ \\
$1.06$ & $19$ & $7$ & $\checkmark$ \\
$1.07$ & $7$  & $6$ & $\checkmark$ \\
$1.08$ & $5$  & $4$ & $\checkmark$ \\
\bottomrule
\end{tabular}
\end{table}

\subsection{Delhi Daily Climate (January 2013 - April 2017)}
\begin{figure}[H]
\centering
\includegraphics[width=1\linewidth]{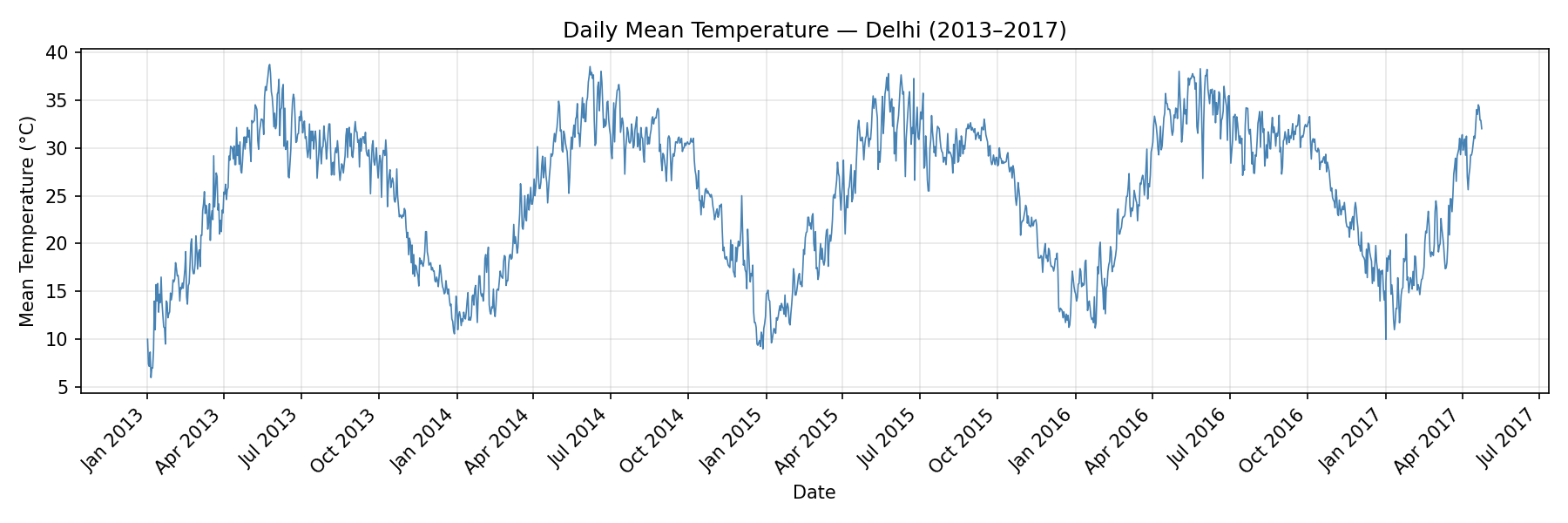}
\caption{Delhi daily mean temperature from January 2013 to April 2017.}
\label{fig:delhi_meantemp}
\end{figure}

We analyze the daily mean temperature of Delhi over $N = 1576$ days (Figure \ref{fig:delhi_meantemp} ). The Fourier ratio of the signal $g$ is
$$\mathrm{FR}(g) = 4.5054,$$ which is much closer to the theoretical minimum of 1 than to the theoretical maximum value of $\sqrt{N} = 39.699$.

\subsubsection{Testing Prediction 1: Small FR Implies Small $\Lambda$}

Consistent with Prediction~1, Table~\ref{tab:climate} shows a small generating set $\Lambda$. Even with $|\Gamma| = 120$ (for $\eta = 1.3$), only $8$ basis elements are needed to span the entire set. This confirms that a small Fourier ratio implies a small $\Lambda$.

\begin{table}[h!]
\centering
\caption{Greedy spanning set $\Lambda$ for the large spectrum $\Gamma$ of the mean temperature time series.}
\begin{tabular}{cccc}
\toprule
$\eta$ & $|\Gamma|$ & $|\Lambda|$ & Spanned 
 \\
\midrule
1.30 & 120 & 8 & \checkmark \\
1.40 & 32  & 8 & \checkmark  \\
1.43 & 14  & 7 & \checkmark  \\
1.45 & 5   & 5 & \checkmark  \\
\bottomrule
\end{tabular}
\label{tab:climate}
\end{table}

\subsubsection{Testing Prediction 3: $\eta$ Scaling}
As $\eta$ increases from $1.3$ to $1.45$, we observe that $|\Lambda|$ decreases from $8$ to $5$. This trend aligns with the theoretical $\eta^{-2}$ scaling.

\subsubsection{Testing Prediction 4: Additive Combination Structure}

Table~\ref{tab:climate} shows that the greedy algorithm successfully spans 
the entire $\Gamma$ in all cases. 
This verifies Prediction~4: every element of $\Gamma$ can be expressed as a 
$\{-1,0,1\}$-linear combination of elements of $\Lambda$.
\begin{figure}[H]
\centering
\includegraphics[width=1\linewidth]{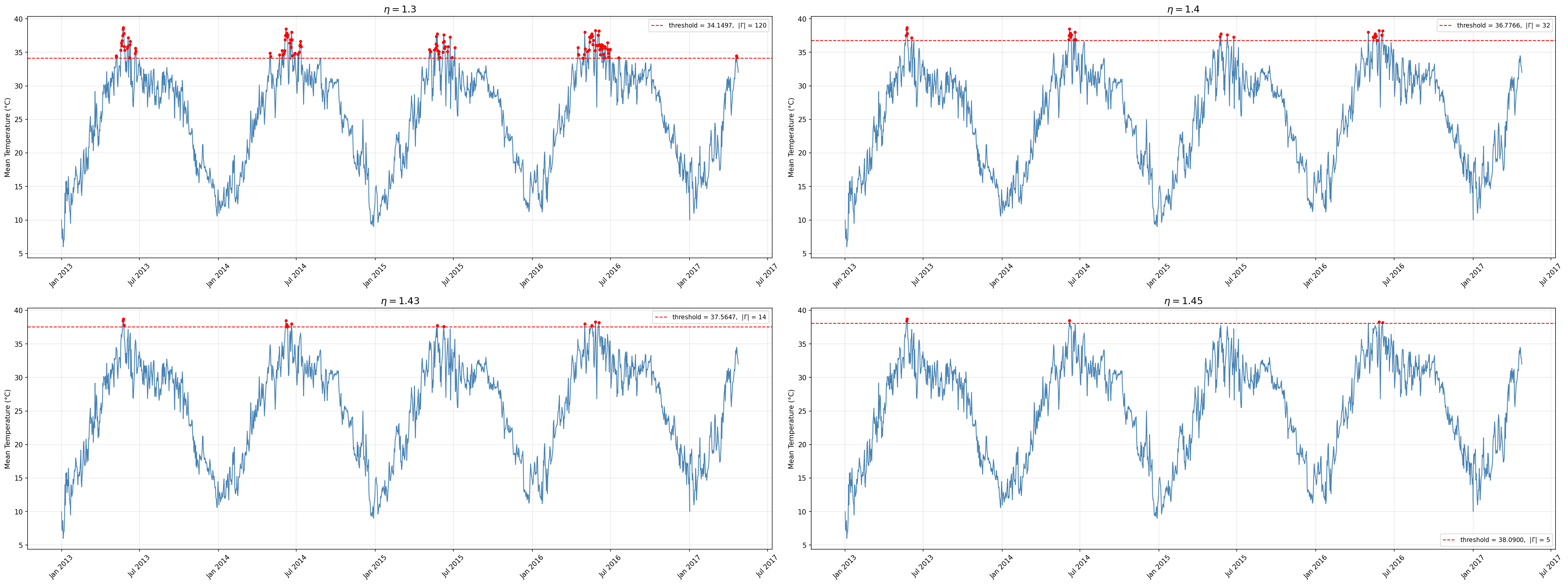}
\caption{Large spectrum $\Gamma$ of the mean temperature time series for $\eta \in \{1.30, 1.40, 1.43, 1.45\}$. Red dots indicate points in $\Gamma$.}
\label{fig:gamma_climate}
\end{figure}

\subsection{Mean-Centered Signals: A Stress Test}
The original signals $f$ (Inflation) and $g$ (Climate) are strictly positive. For such signals, the large spectrum $\Gamma$ captures only the upper tail, which might be trivially structured. To provide a more stringent test, we analyze the mean-centered versions:
\[
f' = f - \overline{f}, \quad  g'=g-\overline{g},
\]
where $\overline{f}$ and $\overline{g}$ are the respective sample means. This creates signals with both positive and negative values. For the inflation data, the Fourier ratio increases dramatically to 
\[
\mathrm{FR}(f') = 10.7853,
\]
which exceeds $\frac{1}{e}\sqrt{N} \approx 8.43$. Thus, the stronger bound 
\eqref{eq:simple_upper_bound} no longer applies, and we must use the looser bound \eqref{eq:general_upper_bound}. In contrast, for the Delhi climate data, the Fourier ratio increases to
$$\mathrm{FR}(g') = 12.6834,$$
but remains just below the value  $\frac{1}{e}\sqrt{1576}\approx 14.60$. Thus, despite the increase in FR, the climate data remains within the regime where the stronger bound \eqref{eq:simple_upper_bound} applies.

\begin{table}[h]
\centering
\caption{Greedy spanning set $\Lambda$ for the large spectrum $\Gamma$ of the mean-centered inflation rate time series.}
\label{tab:inflation_meancent}
\begin{tabular}{cccc}
\toprule
$\eta$ & $|\Gamma|$ & $|\Lambda|$ & Spanned \\
\midrule
$1.5$ & $72$ & $7$ & $\checkmark$ \\
$2.0$ & $25$ & $7$ & $\checkmark$ \\
$2.5$ & $9$  & $6$ & $\checkmark$ \\
$2.7$ & $6$  & $4$ & $\checkmark$ \\
\bottomrule
\end{tabular}
\end{table}
\subsubsection{Testing Prediction 2: Larger FR Regime}

Prediction~2 does not assert that $|\Lambda|$ will increase; rather, it 
acknowledges that the theory provides a weaker guarantee.

Tables~\ref{tab:inflation_meancent} and \ref{tab:delhi_mean_cent} confirm this; the theoretical bounds for both signals become considerably larger than in the raw data cases. However, the empirical cardinality $|\Lambda|$ remains small, ranging from 4 to 7 for inflation data, and 3 to 8 for climate. This stability across both datasets, regardless of whether they satisfy the condition for the stronger bound, suggests that the underlying additive structure is highly robust.

\subsubsection{Testing Prediction 3: Scaling with $\eta$}
For the mean-centered inflation signal, as $\eta$ increases from $1.5$ to 
$2.7$, we observe that $|\Lambda|$ decreases from $7$ to $4$, and for the mean-centered climate signal, as $\eta$ increases from $1.6$ to $2.4$, $|\Lambda|$ decreases from $8$ to $3$, consistent with the predicted trend. 

\subsubsection{Testing Prediction 4: Persistence of Additive Combinations}

\begin{figure}[h]
\centering
\includegraphics[width=1\linewidth]{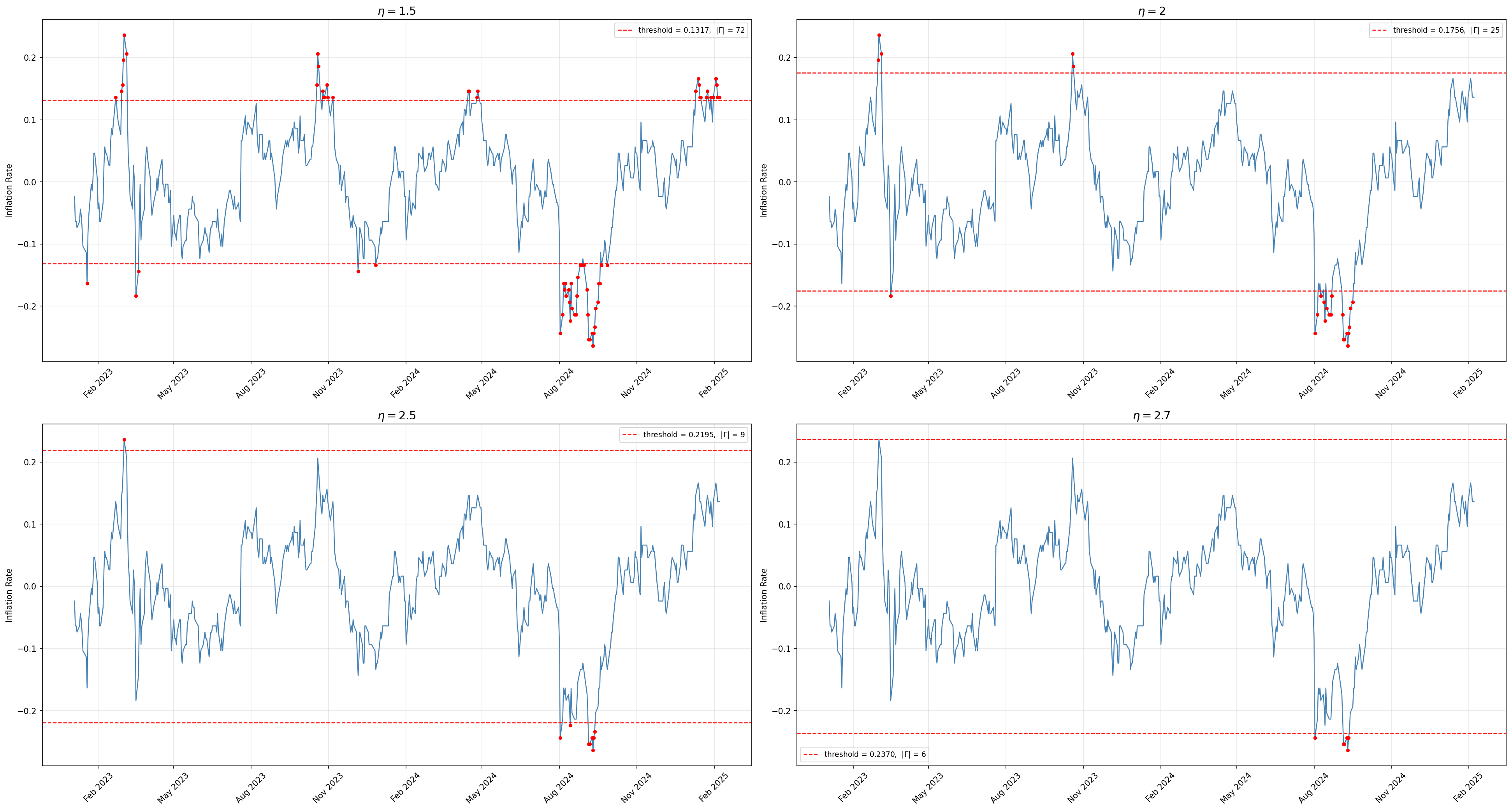}
\caption{Large spectrum $\Gamma$ of the mean-centered inflation rate time series for $\eta \in \{1.5, 2, 2.5, 2.7\}$. Red dots indicate points in $\Gamma$.}
\label{fig:gamma_meancent}
\end{figure}

\begin{figure}[H]
\centering
\includegraphics[width=1\linewidth]{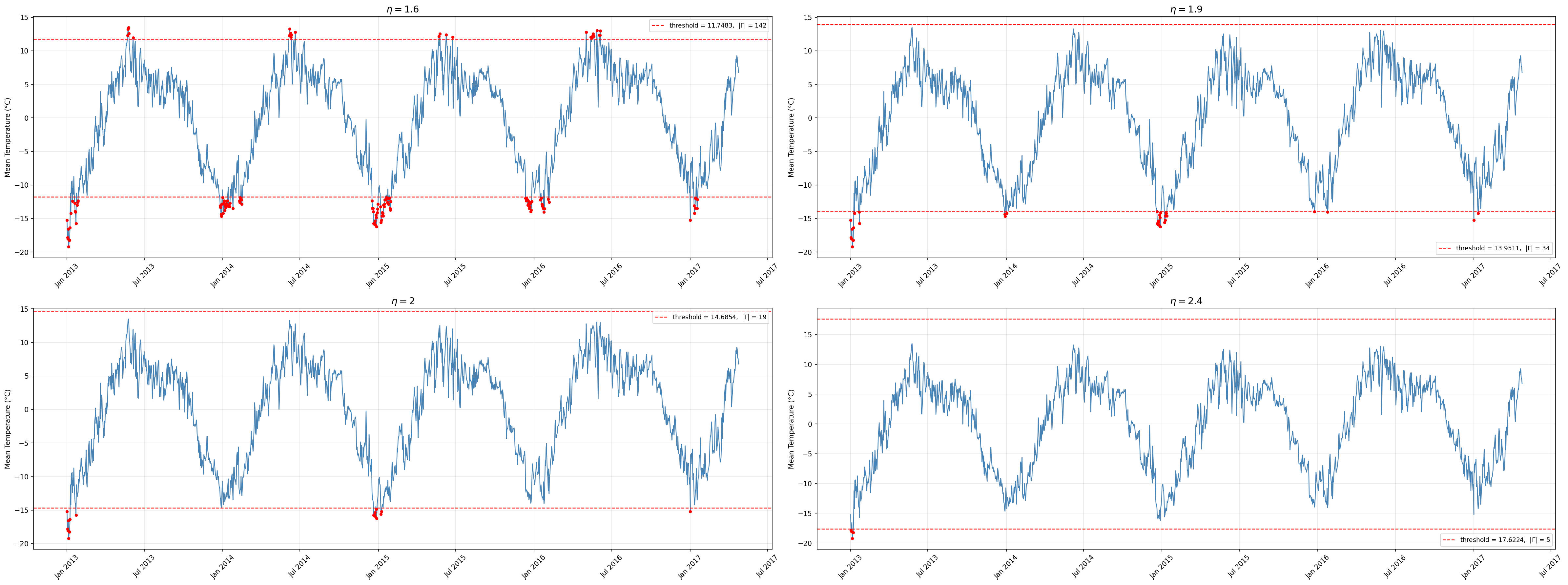}
\caption{Large spectrum $\Gamma$ of the mean-centered mean temperature time series for $\eta \in \{1.6, 1.9, 2, 2.4\}$. Red dots indicate points in $\Gamma$.}
\label{fig:gamma_climate_meancent}
\end{figure}

 The greedy algorithm continues to achieve full coverage of $\Gamma$ in all cases, as seen in the ``Spanned" columns of the tables. Every element of the large spectrum for the mean-centered signals can still be expressed as a $\{-1,0,1\}$-linear combination of the elements in $\Lambda$. While the large spectrum $\Gamma$ becomes visually less clustered in the high-FR regime (Figures~\ref{fig:gamma_meancent} and \ref{fig:gamma_climate_meancent}), the existence of a compact generating set indicates that even when the Fourier ratio is large, the large values of a real-world time series retain significant additive 
structure.

\begin{table}
\centering
\caption{Greedy spanning set $\Lambda$ for the large spectrum $\Gamma$ of the mean-centered mean temperature time series.}
\begin{tabular}{cccc}
\toprule
$\eta$ & $|\Gamma|$ & $|\Lambda|$ & Spanned  \\
\midrule
1.60 & 142 & 8 & \checkmark  \\
1.90 & 34 & 7 & \checkmark \\
2.00 & 19  & 6 & \checkmark  \\
2.40 & 5   & 3 & \checkmark   \\
\bottomrule
\end{tabular}
\label{tab:delhi_mean_cent}
\end{table}

\subsection{Summary of Validation}
\label{subsec:validation}

Table~\ref{tab:validation} summarizes how each theoretical prediction fared 
against the empirical data.

\begin{table}[H]
\centering
\caption{Validation of theoretical predictions.}
\label{tab:validation}
\begin{tabular}{p{4cm}p{4cm}p{4cm}c}
\toprule
\textbf{Prediction} & \textbf{Inflation ($N=526$)} & \textbf{Climate ($N=1576$)} & \textbf{Status} \\
\midrule
\textbf{P1:} Small $\mathrm{FR} \implies$ small $|\Lambda|$ & $|\Lambda| \in [4, 7]$ for $\mathrm{FR} = 1.41$ & $|\Lambda| \in [5, 8]$ for $\mathrm{FR} = 4.51$ & Confirmed \\
\addlinespace[0.5em]
\textbf{P2:} High $\mathrm{FR} \implies$ looser bound & Bound increases, but $|\Lambda|$ remains $\leq 7$ & Bound increases, but $|\Lambda|$ remains $\leq 8$ & Confirmed \\
\addlinespace[0.5em]
\textbf{P3:} $|\Lambda|$ scales with $\eta$ & $|\Lambda|$ decreases from $7$ to $4$ as $\eta$ increases & $|\Lambda|$ decreases from $8$ to $3$ as $\eta$ increases & Confirmed \\
\addlinespace[0.5em]
\textbf{P4:} Additive structure & greedy algorithm achieved full spanning & greedy algorithm achieved full spanning & Confirmed \\
\bottomrule
\end{tabular}
\end{table}

\subsection{Discussion}
\label{subsec:discussion}

The numerical experiments above validate the four theoretical predictions 
derived from Corollary~\ref{cor:mama}. The key findings are:

\begin{enumerate}
\item \textbf{Small generating sets:} In all cases tested, the size of $\Lambda$ remained an order of magnitude smaller than $|\Gamma|$, often as small as $4$--$7$ even when $|\Gamma|$ contained dozens of points.

\item \textbf{Full spanning:} The greedy algorithm successfully spanned the entire large spectrum $\Gamma$ in every instance, indicating that the additive structure guaranteed by Chang's Lemma is not merely existential but can be efficiently realized.

\item \textbf{Robustness to mean-centering:} Even after removing the mean---which dramatically increased the Fourier ratio and moved us into a regime with weaker theoretical guarantees---the additive structure persisted. This suggests that the phenomenon may be even more general than our current bounds capture.

\item \textbf{$\eta$ scaling:} The observed $|\Lambda|$ decreased with increasing $\eta$ as predicted by the $\eta^{-2}$ scaling in the theoretical bound.
\end{enumerate}

These results provide empirical validation for the theoretical framework developed in Section~\ref{section:introduction}. The large values of real-world time series are not merely outliers but exhibit genuine additive structure that can be exploited for compression and analysis. The Fourier ratio serves as a practical diagnostic: when $\mathrm{FR}(f)$ is small, the theory guarantees strong additive structure; when it is larger, the structure may still be present, but our current bounds become less informative.

\subsection{Limitations and Caveats}
\label{subsec:limitations}

Several limitations of the present study should be acknowledged.

\paragraph{Domain coverage.}
Our numerical experiments focus on US inflation and Delhi daily climate data. While these represent two distinct and physically different domains (economic and meteorological), the conclusions may not generalize to all classes of time series (e.g., financial returns with different tail behavior, physiological signals, or high-frequency trading data). Testing on a broader corpus of time series is an important direction for future work.

\paragraph{Unknown constant $C$.}
The theoretical bound in Corollary~\ref{cor:mama} contains an unspecified absolute constant $C$. Our numerical validation shows that $|\Lambda|$ is smaller than $\mathrm{Bound}/C$ for all tested parameters, which is consistent with the existence of such a $C$, but we cannot determine its value or tightness from these experiments alone.

\paragraph{Greedy algorithm non-optimality.}
Algorithm~\ref{alg:greedy} produces some $\Lambda$ that spans $\Gamma$, but not necessarily the minimal such set. The true minimal $|\Lambda|$ may be smaller than reported. Finding the minimal $\Lambda$ is a set cover problem that is computationally hard in general, though our greedy heuristic suffices for validation.

\paragraph{Mean-centering as a stress test.}
The mean-centered signals $f' = f - \overline{f}$ and $g' = g - \overline{g}$ have a much larger Fourier ratio, but this transformation is not typical in practice. Our purpose was to test the theory in a regime where its guarantees weaken, not to recommend mean-centering as a preprocessing step. Practitioners should compute $\mathrm{FR}(f)$ on their original data to determine whether the strong bound applies.

\section{Conclusion and Future Work}
\label{sec:conclusion}

We have shown that additive combinatorics, specifically the generalized 
Chang's lemma and the Fourier ratio, provides a rigorous foundation for the 
empirical observation that large values of time series exhibit meaningful 
structure. The main theoretical result (Corollary~\ref{cor:mama}) establishes 
that when $\mathrm{FR}(f)$ is small, the set of large values is additively 
generated by a small set $\Lambda$. 

Our numerical experiments on US inflation data and Delhi daily mean temperature records validated the four predictions derived from this theory:
\begin{itemize}

\item For the original signals (FR=1.41 for inflation; FR=4.51 for climate), we found that $|\Lambda|$ remained between 4 and 8, even when  $|\Gamma|$ contained over 100 points.
\item After mean-centering ($\mathrm{FR}=10.79$ and $12.68$, respectively), the additive structure persisted. This occurred even when the inflation data crossed the critical value for the stronger theoretical bound, suggesting that the phenomenon is more robust than the current theory guarantees.
\item The observed $|\Lambda|$ decreased with $\eta$ as predicted.
\item Full $\{-1,0,1\}$-spanning was achieved in all cases.
\end{itemize}

These results have several implications for time series analysis. First, they 
provide a justification for why preprocessing techniques like Box--Cox 
transforms (\cite{BoxCox1964}) work: they reduce the Fourier ratio, making the additive structure 
more accessible. Second, they suggest a new compression strategy: store only 
the generating set $\Lambda$ and reconstruct the large values via additive 
combinations. Third, they offer a diagnostic: when $\mathrm{FR}(f)$ is small, 
one can confidently apply methods that exploit tail structure.

\vskip.125in 

Future work includes:
\begin{itemize}
\item Testing the predictions on additional datasets (retail sales, 
financial returns, network traffic) to assess the possible nuanced differences between various industries. 
\item Developing algorithms that find the \emph{minimal} $\Lambda$ (our greedy 
algorithm provides only an upper bound).
\item Exploring whether the additive structure can be exploited for forecasting, 
not just compression.
\item Tightening the bound in the large-FR regime, where our experiments suggest the theory is currently loose.
\item Investigating connections to other notions of structure from additive combinatorics, such as Freiman's theorem (\cite{Freiman1966}) and sumset estimates.
\end{itemize}

Finally, this paper opens a broader question: what other notions of 
``structure'' from additive combinatorics can be brought to bear on time 
series analysis? We believe this cross-pollination has significant untapped 
potential.

\subsection*{Competing Interests}

The authors declare no competing interests.

\end{document}